\documentclass[11pt,a4paper]{article}
\usepackage[margin=1in]{geometry}
\usepackage{latexsym, amsfonts, amsmath,tikz}
\usepgflibrary{arrows}

\newcommand{\be}{\begin{equation}}
\newcommand{\ee}{\end{equation}}
\newcommand{\bea}{\begin{eqnarray}}
\newcommand{\eea}{\end{eqnarray}}
\newcommand{\bean}{\begin{eqnarray*}} 
\newcommand{\eean}{\end{eqnarray*}}

\newcommand{\brray}{\begin{array}}
\newcommand{\erray}{\end{array}}
\newcommand{\ben}{\begin{equation}{nonumber}}
\newcommand{\een}{\end{equation}{nonumber}}

\newtheorem{dfn}{Definition}[section]
\newtheorem{thm}[dfn]{Theorem}
\newtheorem{lmma}[dfn]{Lemma}
\newtheorem{ppsn}[dfn]{Proposition}
\newtheorem{crlre}[dfn]{Corollary}
\newtheorem{xmpl}[dfn]{Example}
\newtheorem{rmrk}[dfn]{Remark}
\newcommand{\bdfn}{\begin{dfn}}
\newcommand{\bthm}{\begin{thm}}
\newcommand{\blmma}{\begin{lmma}}
\newcommand{\bppsn}{\begin{ppsn}}
\newcommand{\bcrlre}{\begin{crlre}}
\newcommand{\bxmpl}{\begin{xmpl}}
\newcommand{\brmrk}{\begin{rmrk}}
\newcommand{\edfn}{\end{dfn}}
\newcommand{\ethm}{\end{thm}}
\newcommand{\elmma}{\end{lmma}}
\newcommand{\eppsn}{\end{ppsn}}
\newcommand{\ecrlre}{\end{crlre}}
\newcommand{\exmpl}{\end{xmpl}}
\newcommand{\ermrk}{\end{rmrk}}

\newcommand{\IC}{\mathbb{C}}

\newcommand{\IR}{\mathbb{R}}

\newcommand{\cla}{{\cal A}}
\newcommand{\clb}{{\cal B}}
\newcommand{\clc}{{\cal C}}
\newcommand{\cld}{{\cal D}}
\newcommand{\cle}{{\cal E}}
\newcommand{\clf}{{\cal F}}

\newcommand{\clh}{{\cal H}}
\newcommand{\cli}{{\cal I}}

\newcommand{\cll}{{\cal L}}
\newcommand{\clm}{{\cal M}}
\newcommand{\cln}{{\cal N}}

\newcommand{\clp}{{\cal P}}
\newcommand{\clq}{{\cal Q}}

\newcommand{\cls}{{\cal S}}

\newcommand{\clu}{{\cal U}}
\newcommand{\clv}{{\cal V}}

\def\a*{{\cal A}_{h,*}}
\def\B{{\cal B}(h)}
\def\B1{{\cal B}_1(h)}
\def\b{{\cal B}^{\rm s.a.}(h)}
\def\b1{{\cal B}^{\rm s.a.}_1(h)}

\newcommand{\ot}{\otimes}

\newcommand{\raro}{\rightarrow}

\newcommand{\lgl}{\langle}
\newcommand{\rgl}{\rangle}

\newcommand{\midarrow}{\tikz \draw[-triangle 90] (0,0) -- +(.1,0);}
\def \qed {$\Box$}

\def\a*{{\cal A}_{h,*}}
\def\B{{\cal B}(h)}
\def\B1{{\cal B}_1(h)}
\def\b{{\cal B}^{\rm s.a.}(h)}
\def\b1{{\cal B}^{\rm s.a.}_1(h)}

\begin{document}
\begin{center}
{\Large{\bf Quantum Symmetry of Graph $C^{\ast}$-algebras associated with connected Graphs}}
\vspace{0.1in}\\
{\large {\bf Soumalya Joardar \footnote{Acknowledges support from INSPIRE Faculty grant funded by Department of Science and Technology, India}} and \bf Arnab Mandal \footnote{Acknowledges SERB (Govt.of India)}}
\end{center}
\begin{abstract}
We define a notion of quantum automorphism groups of graph $C^{\ast}$-algebras for finite, connected graphs. Under the assumption that the underlying graph does not have any multiple edge or loop, the quantum automorphism group of the underlying directed graph in the sense of T. Banica (\cite{Ban}) (which is also the symmetry object in the sense of \cite{Web}) is shown to be a quantum subgroup of quantum automorphism group in our sense. Quantum symmetries for some concrete graph $C^{\ast}$-algebras have been computed. 
\end{abstract}
\section{Introduction}

\indent With the advent of noncommutative geometry it was natural for Mathematicians to find the right notion of symmetry for `noncommutative' spaces. In fact in 1995, at Les Houches summer school on quantum symmetries, Alain Connes raised the question of finding a proper notion of quantum automorphism group of classical or non commutative spaces. It was clear that such symmetry had to be a generalization of classical group symmetries for classical spaces. In 1998, S. Wang came up with a seminal work on this problem in \cite{Wang}, where he defined the notion of quantum automorphism group of finite spaces. The idea was to formulate the classical group symmetry problem in categorical language i.e. to realize the automorphism group of the underlying space as a universal object in a certain category of compact groups and then letting the objects in the category be replaced by compact quantum groups and the group homomorphisms by CQG morphisms. Then the main challenge was to prove the existence of universal object in the larger category. As it turned out even for an $n$-point space where the underlying $C^{\ast}$-algebra is a finite dimensional $C^{\ast}$-algebra, the quantum automorphism is significantly larger than the classical symmetry group $S_{n}$ and in fact as a $C^{\ast}$-algebra the quantum automorphism group is infinite dimensional for $n\geq 4$. In that paper he also defined quantum automorphism group for finite dimensional $C^{\ast}$-algebras. And in deed for such $C^{\ast}$-algebras, the universal object fails to exist in general. S. Wang remedied this problem by his pioneering idea of restricting the category of $C^{\ast}$-actions of CQG's by demanding that such an action should also preserve some suitable linear functional. Since that paper, study of quantum symmetries has come a long way. Quantum symmetry of finite structures was extended for finite graphs by T. Banica, J. Bichon and others (\cite{Ban},\cite{Bichon} and references therein). Later on D. Goswami and others defined quantum symmetries for continuous structures in the framework of noncommutative geometry (\cite{Debashish},\cite{Laplace} and references therein).\\
\indent It is interesting to see that even if for a directed graph the function algebra over the graph is finite dimensional $C^{\ast}$-algebra, the associated graph $C^{\ast}$-algebra is not so (see \cite{Pask}) in general. So in the context of quantum symmetry it is interesting to study the quantum symmetries of the graph $C^{\ast}$-algebras. It is expected that we shall have larger quantum symmetry than that of the underlying graph. In this paper we define a notion of quantum symmetry of graph $C^{\ast}$-algebras. We take the same approach as S. Wang and define a category which is more restrictive than a CQG only having a faithful action on the graph $C^{\ast}$-algebra (much like the case of finite dimensional $C^{\ast}$-algebra). Then we show for finite, connected graphs, universal object exists in the category and call them quantum symmetry of the graph $C^{\ast}$-algebra. As expected we show that the quantum symmetry is quite large and if we also assume that the underlying graph has no multiple edge or loop, the quantum automorphism group of the underlying directed graph in the sense of T. Banica always is an object in the category. In general the quantum symmetry group of the graph $C^{\ast}$-algebra is strictly larger than the quantum symmetry group of underlying graph in the sense of Banica. We show it by computing the quantum symmetry group of a complete graph with 2 vertices. Note that recently M. Weber and S. Schmidt also studied the quantum symmetries of graph $C^{\ast}$-algebras (see \cite{Web}). They have taken a category which is more of an algebraic category where as our category is more akin to isometric formulation. Also in \cite{Ortho}, A. Skalski and T. Banica defined a notion of quantum symmetry for Cuntz-Krieger algebras. They have taken the category of CQG actions preserving a given orthogonal filtration corresponding to faithful state on Cuntz-Krieger algebra.\\
\indent Now we briefly discuss the organization of the paper. We start with a preliminary section where we collect some standard facts about compact quantum groups and quantum automorphism groups. Moreover, we also recall the quantum symmetry of directed graphs in the sense of T. Banica from \cite{Ban} and \cite{Fulton}. 4th section contains the notion of quantum automorphism group of graph $C^{\ast}$-algebra and  a non rigidity result. On the way we also show that in general the quantum symmetry group of a directed and connected graph without any multiple edge or loop in the sense of T. Banica is an object in our category. In the last section, the quantum automorphism groups of a few finite, connected graphs are computed.  
\section{Notations and conventions}
 All the $C^{\ast}$-algebras are unital. For a $C^{\ast}$-algebra $\cla$, $\cla^{\ast}$ will denote the space of all bounded linear functionals on $\cla$. The tensor product $\ot$ between two $C^{\ast}$-algebras are injective tensor products. ${\rm id}$ will denote the identity map. For a set $A$, ${\rm Sp}(A) \ ({\rm resp.} \ \overline{\rm Sp}(A))$ will denote the linear span (resp. closed linear span) of the elements of $A$. By a universal object in a category we shall always mean the initial object in the category.
\section{Preliminaries}
\subsection{Compact quantum groups and quantum automorphism groups}
\label{qaut}
In this subsection we recall the basics of compact quantum groups and their actions on $C^{\ast}$-algebras. The facts collected in this Subsection are well known and we refer the readers to \cite{Van}, \cite{Woro}, \cite{Wang} for details. 
\bdfn
A compact quantum group (CQG) is a pair $(\clq,\Delta)$ such that $\clq$ is a unital $C^{\ast}$-algebra and $\Delta:\clq\raro\clq\ot\clq$ is a unital $C^{\ast}$-homomorphism satisfying\\
(i) $({\rm id}\ot\Delta)\Delta=(\Delta\ot{\rm id})\Delta$.\\
(ii) {\rm Sp}$\{\Delta(\clq)(1\ot\clq)\}$ and {\rm Sp}$\{\Delta(\clq)(\clq\ot 1)\}$ are dense in $\clq\ot\clq$.
\edfn
Given a CQG $\clq$, there is a canonical dense Hopf $\ast$-algebra $\clq_{0}$ in $\clq$ on which an antipode $\kappa$ and counit $\epsilon$ are defined. Given two CQG's $(\clq_{1},\Delta_{1})$ and $(\clq_{2},\Delta_{2})$, a CQG morphism between them is a $C^{\ast}$-homomorphism $\pi:\clq_{1}\raro\clq_{2}$ such that $(\pi\ot\pi)\Delta_{1}=\Delta_{2}\pi$.
\vspace{0.1in}\\
{\it Examples}:\\
1. Let $U_{n}^{+}$ be the universal $C^{\ast}$-algebra generated by $n^2$ elements $\{q_{ij}\}_{i,j=1,\ldots,n}$ such that both the matrices $((q_{ij}))$ and $((q_{ij}))^t$ are unitaries. Then $U_{n}^{+}$ is a CQG with coproduct given on the generators by $\Delta_{0}(q_{ij})=\sum_{k=1}^{n}q_{ik}\ot q_{kj}$ (see \cite{Wangfree} where this CQG has been denoted by $A_{u}(n)$).
\vspace{0.1in}\\
2. $(C(S^{1}),\Delta)$ is a commutative CQG with coproduct given on the generator $z$ by diagonal action i.e. $\Delta(z)=z\ot z$. The free product of $n$ copies of $C(S^{1})$ is a non commutative CQG (see \cite{Wangfree} for generalities on free product of CQG's) denoted by $(\underbrace{C(S^{1})\ast\ldots\ast C(S^{1})}_{n-times},\Delta)$ where the coproduct on the generators is given by the diagonal action i.e. $\Delta(z_{i})=z_{i}\ot z_{i}$ for all $i=1,\ldots,n$. This free product has the following universal property:\\
If $(\clq,\Delta^{\prime})$ is a CQG generated by $n$ unitaries $q_{i}$ such that $\Delta^{\prime}(q_{i})=q_{i}\ot q_{i}$ for all $i$, then there is a surjective CQG morphism from $(\underbrace{C(S^{1})\ast\ldots\ast C(S^{1})}_{n-times},\Delta)$ onto $(\clq,\Delta^{\prime})$ sending $z_{i}$ to $q_{i}$.
\vspace{0.1in}\\
3. For an invertible $n\times n$ matrix $Q\in M_{n}(\mathbb{C})$, let $A_{u^{t}}(Q)$ be the universal $C^{\ast}$-algebra generated by elements $(u_{ij})_{i,j=1,\ldots,n}$ satisfying the following relations
\begin{eqnarray*}
 U^{t^{\ast}}U^{t}=U^{t}U^{t^{\ast}}={\rm Id}_{n\times n}, \ UQ^{-1}U^{\ast}Q=Q^{-1}U^{\ast}QU={\rm Id}_{n\times n},
\end{eqnarray*}
where $U$ is the $A_{u^{t}}(Q)$-valued $n\times n$ matrix $((u_{ij}))$. Then adapting the arguments of Theorem 1.3 of \cite{Wang_univ} it can be shown that $A_{u^{t}}(Q)$ has a CQG structure with coproduct given on the generators by $\Delta_{0}(u_{ij})=\sum_{k=1}^{n}u_{ik}\ot u_{kj}$. Note that for $Q={\rm Id}_{n\times n}$, $A_{u^{t}}(Q)$ is $U^{+}_{n}$. Also it is easy to see that as a CQG, $(A_{u^{t}}(Q),\Delta_{0})$ is isomorphic to $(A_{u}(Q^{-1}),\Delta^{\rm cop})$ (For the definition of $(A_{u}(Q^{-1}),\Delta)$ see \cite{Wang_univ}). 
\vspace{0.1in}\\
\indent Before going to faithful $C^{\ast}$-action of CQG's we briefly recall the doubling procedure of a compact quantum group from \cite{Soltan} for future use in this paper. Let $(\mathcal{Q},\Delta)$ be a CQG with a CQG-automorphism $\theta$ such that $\theta^2={ \rm Id}$. The doubling of this CQG, 
 say $(\mathcal{D}_{\theta}(\mathcal{Q}),\tilde{\Delta})$, is  given by  $\mathcal{D}_{\theta}(\mathcal{Q}) :=\mathcal{Q}\oplus \mathcal{Q}$ (direct sum as a $C^{\ast}$-algebra),
 and the coproduct is defined by the following
 \begin{eqnarray*}
&& \tilde{\Delta} \circ \xi= (\xi \ot \xi + \eta \ot [\eta \circ \theta])\circ  \Delta\\
&& \tilde{\Delta} \circ \eta= (\xi \ot \eta + \eta \ot [\xi \circ \theta])\circ \Delta,
\end{eqnarray*}
 where we have denoted  the injections of $\mathcal{Q}$ onto 
 the first and second coordinate in $\mathcal{D}_{\theta}(\mathcal{Q})$ by $\xi$ and $\eta$ respectively, i.e. 
$\xi(a)=(a,0), \  \eta(a)= (0,a), \ (a \in \mathcal{Q}).$
\bdfn
Given a (unital) $C^{\ast}$-algebra $\clc$, a CQG $(\clq,\Delta)$ is said to act faithfully on $\clc$ if there is a unital $C^{\ast}$-homomorphism $\alpha:\clc\raro\clc\ot\clq$ satisfying\\
(i) $(\alpha\ot {\rm id})\alpha=({\rm id}\ot \Delta)\alpha$.\\
(ii) {\rm Sp}$\{\alpha(\clc)(1\ot\clq)\}$ is dense in $\clc\ot\clq$.\\
(iii) The $\ast$-algebra generated by the set  $\{(\omega\ot{\rm id})\alpha(\clc): \omega\in\clc^{\ast}\}$ is norm-dense in $\clq$.
\edfn
For a faithful action of a CQG $(\clq,\Delta)$ on a unital $C^{\ast}$-algebra $\clc$, there is a norm dense $\ast$-subalgebra $\clc_{0}$ of $\clc$ such that the canonical Hopf-algebra $\clq_{0}$ coacts on $\clc_{0}$.
\bdfn
(Def 2.1 of \cite{Bichon})
Given a unital $C^{\ast}$-algebra $\clc$, quantum automorphism group of $\clc$ is a CQG $(\clq,\Delta)$ acting faithfully on $\clc$ satisfying the following universal property:\\
\indent If $(\clb,\Delta_{\clb})$ is any CQG acting faithfully on $\clc$, there is a surjective CQG morphism $\pi:\clq\raro\clb$ such that $({\rm id}\ot \pi)\alpha=\beta$, where $\beta:\clc\raro\clc\ot\clb$ is the corresponding action of $(\clb,\Delta_{\clb})$ on $\clc$ and $\alpha$ is the action of $(\clq,\Delta)$ on $\clc$.
\edfn
\brmrk
In general universal object might fail to exist in the above category. For the existence of universal object one generally restricts the category. In addition to $\alpha$ being a faithful action it is common to add some kind of `volume' preserving condition i.e. one assumes some linear functional $\tau$ on $\clc$ such that $(\tau\ot {\rm id})\alpha(a)=\tau(a).1$ for $a$ in some suitable subspace of $\clc$ (see \cite{Wang}, \cite{Debashish}).
\ermrk
{\it Example}:\\
1. If we take a space of $n$ points $X_{n}$ then the quantum automorphism group of the $C^{\ast}$-algebra $C(X_{n})$ is given by the CQG (denoted by $S_{n}^{+}$) which as a $C^{\ast}$-algebra is the universal $C^{\ast}$ algebra generated by $\{u_{ij}\}_{i,j=1,\ldots,n}$ satisfying the following relations (see Theorem 3.1 of \cite{Wang}):
\begin{displaymath}
 u_{ij}^{2}=u_{ij},u_{ij}^{\ast}=u_{ij},\sum_{k=1}^{n}u_{ik}=\sum_{k=1}^{n}u_{kj}=1, \ i,j=1,\ldots,n.
\end{displaymath}
The coproduct on the generators is given by $\Delta(u_{ij})=\sum_{k=1}^{n}u_{ik}\ot u_{kj}$.\vspace{0.1in}\\
2. If we take the $C^{\ast}$-algebra $M_{n}(\mathbb{C})$, then as remarked earlier if we take the category of CQG's only having faithful $C^{\ast}$-action on $M_{n}(\mathbb{C})$, then universal object fails to exist in the category. But in addition if we assume any object in the category also has to preserve a linear functional $\phi$ on $M_{n}(\mathbb{C})$, then existence of the universal object can be shown (see \cite{Wang}). 
\subsection{Quantum symmetry of finite graphs without loops and multiple edges}
\label{graph_bich}
A {\bf finite} directed graph is a collection of finitely many edges and vertices. If we denote the edge set of a graph $\Gamma$ by $E=(e_{1},\ldots,e_{n})$ and set of vertices of $\Gamma$ by $V=(v_{1},\ldots,v_{m})$ then recall the maps $s,t:E\raro V$ and the adjacency matrix $D$ whose $ij$-th entry is $1$ if the ordered pair $(i,j)\in E$ and $0$ otherwise (see page $8$ of \cite{Fulton}).
\bdfn
$\Gamma$ is said to be without multiple edges if $s(e_{i})=s(e_{j})$ and $t(e_{i})=t(e_{j})$ imply that $i=j$. $\Gamma$ has no loop if $s(e_{i})\neq t(e_{i})$ for all $i=1,\ldots,n$.
\edfn
\bdfn
$(Q^{aut}_{Ban}(\Gamma),\Delta)$ for a graph $\Gamma$ without multiple edge or loop is defined to be the quotient $S^{+}_{n}/(AD-DA)$, where $A=((u_{ij}))_{i.j=1,\ldots,m}$, and $D$ is the adjacency matrix for $\Gamma$. The coproduct on the generators is again given by $\Delta(u_{ij})=\sum_{k=1}^{m}u_{ik}\ot u_{kj}$.
\edfn
\bthm
(Lemma 3.1.1 of \cite{Fulton}) The quantum automorphism group $(Q^{aut}_{Ban}(\Gamma),\Delta)$ of a finite graph $\Gamma$ with $n$ edges and $m$ vertices (without loop or multiple edge) is the universal $C^{\ast}$-algebra generated by $(u_{ij})_{i,j=1,\ldots,m}$ satisfying the following relations:
\begin{eqnarray}
 && u_{ij}^{\ast}=u_{ij}, u_{ij}u_{ik}=\delta_{jk}u_{ij}, u_{ji}u_{ki}=\delta_{jk}u_{ji}, \sum_{l=1}^{m}u_{il}=\sum_{l=1}^{m}u_{li}=1, 1\leq i,j,k\leq m \label{1}\\
 && u_{s(e_{j})i}u_{t(e_{j})k}=u_{t(e_{j})k}u_{s(e_{j})i}=u_{is(e_{j})}u_{kt(e_{j})}=u_{kt(e_{j})}u_{is(e_{j})}=0, e_{j}\in E, (i,k)\not\in E \label{2}
\end{eqnarray}
where the coproduct on the generators is given by $\Delta(u_{ij})=\sum_{k=1}^{m}u_{ik}\ot u_{kj}$.
\ethm
We recall the following Lemma from \cite{Web} (Lemma 3.7):
\blmma
\label{Weber}
Let $\Gamma=(E,V)$ be a finite graph without multiple edges or loops and $e_{j}\in E$. Then for any $v\in V$ such that $s^{-1}(v)=\emptyset$, $u_{vs(e_{j})}=u_{s(e_{j})v}=0$ for all $j=1,\ldots,n$, where $u_{ij}$'s are generators of $(Q^{aut}_{Ban}(\Gamma),\Delta)$. 
\elmma

\section{Quantum symmetry of Graph $C^{\ast}$-algebras}
\subsection{Graph $C^{\ast}$-algebras}
\label{graph}
In this subsection we recall few basic facts of graph $C^{\ast}$-algebras. Readers are referred to \cite{Raeburn}, \cite{Pask} for details. Let $\Gamma=\{E=(e_{1},...,e_{n}),V=(v_{1},...,v_{m})\}$ be a finite graph. A graph is said to be {\bf connected} if every vertex is a source or a target of an edge. In other words for all $v_{l}\in V$, $s^{-1}(v_{l})$ or $t^{-1}(v_{l})$ is non empty. In this paper all the graphs are {\bf finite, connected}.
\bdfn
\label{Graph}
The graph $C^{\ast}$-algebra $C^{\ast}(\Gamma)$ is defined as the universal $C^{\ast}$-algebra generated by partial isometries $\{S_{e_{i}}\}_{i=1,\ldots,n}$ and mutually orthogonal projections $\{p_{v_{k}}\}_{k=1,\ldots,m}$ satisfying the following relations:
\begin{displaymath}
  S_{e_{i}}^{\ast}S_{e_{i}}=p_{t(e_{i})}, \sum_{s(e_{j})=v_{l}}S_{e_{j}}S_{e_{j}}^{\ast}=p_{v_{l}}.
	\end{displaymath}
\edfn
In a graph $C^{\ast}$-algebra $C^{\ast}(\Gamma)$, we have the following (see Subsection 2.1 of \cite{Pask}):\\
1. $\sum_{k=1}^{m}p_{v_{k}}=1$ and $S_{e_{i}}^{\ast}S_{e_{j}}=0$ for $i\neq j$.
\vspace{0.05in}\\
2. $S_{\mu}=S_{e_{1}}S_{e_{2}}\ldots S_{e_{l}}$ is non zero if and only if $\mu=e_{1}e_{2}\ldots e_{l}$ is a path i.e. $t(e_{i})=s(e_{i+1})$ for $i=1,\ldots,(l-1)$.
\vspace{0.05in}\\
3. $C^{\ast}(\Gamma)={\overline{\rm Sp}}\{S_{\mu}S_{\nu}^{\ast}:t(\mu)=t(\nu)\}$.\\
\indent  Let $\cli\subset\{1,...,m\}$ be such that for each $k\in\cli$, $v_{k}$ is not a source of any edge of the graph. Also let $E^{\prime}\subset E\times E$ be such that $S_{e_{i}}S_{e_{j}}^{\ast}\neq 0$ for $(e_{i},e_{j})\in E^{\prime}$. Note that $(e_{i},e_{i})\in E^{\prime}$ for all $i=1,\ldots,n$. Denote the vector space  ${\rm Sp}\{S_{e_{i}}S_{e_{j}}^{\ast},(p_{v_{k}})_{k\in\cli},i,j=1,\ldots,n\}$ by $\clv_{2,+}$. Then we have the following
\blmma
The set $\clb=\{S_{e_{i}}S_{e_{j}}^{\ast},p_{v_{k}}: (e_{i},e_{j})\in E^{\prime}, k\in\cli\}$ is a basis for the vector space $\clv_{2,+}$.
\elmma
{\bf Proof}:\\
Clearly $\clb$ spans $\clv_{2,+}$. So we need to show the linear independence of the set to finish the proof of the Lemma. Let $\Lambda$ be the set of pairs $(i,j)$ such that $(e_{i},e_{j})\in E^{\prime}$. Let $\sum_{(i,j)\in\Lambda}\lambda_{ij}S_{e_{i}}S_{e_{j}}^{\ast}+\sum_{k\in\cli}\mu_{k}p_{v_{k}}=0$. For $l\in\cli$, say there is an $e_{m}$ such that $t(e_{m})=v_{l}$ i.e. $S_{e_{m}}^{\ast}S_{e_{m}}=p_{v_{l}}$. Multiplying the above equation by $p_{v_{l}}$ from right, we get
\begin{displaymath}
\mu_{l}p_{v_{l}}+\sum_{(i,j)\in\Lambda}\lambda_{ij}S_{e_{i}}(S_{e_{m}}S_{e_{j}})^{\ast}S_{e_{m}}=0.
\end{displaymath}
 $S_{e_{m}}S_{e_{j}}$ is non zero if and only if $v_{l}=t(e_{m})=s(e_{j})$ (see after Definition (\ref{Graph})). But by our assumption, $v_{l}$ is not a source of any edge implying that $S_{e_{m}}S_{e_{j}}=0$ for all $j$. Hence we get $\mu_{l}=0$ for all $l\in\cli$. So we are left to prove that $\lambda_{ij}=0$ for all $(i,j)\in\Lambda$. To that end let $(k,l)\in\Lambda$. Multiplying the equation $\sum_{(i,j)\in\Lambda}\lambda_{ij}S_{e_{i}}S_{e_{j}}^{\ast}=0$ by $S_{e_{k}}S_{e_{k}}^{\ast}$ from left and $S_{e_{l}}S_{e_{l}}^{\ast}$ from right and using $S_{e_{i}}^{\ast}S_{e_{j}}=0$ for $i\neq j$, we get $\lambda_{kl}S_{e_{k}}S_{e_{l}}^{\ast}=0.$ Since $S_{e_{k}}S_{e_{l}}^{\ast}\neq 0$, we get $\lambda_{kl}=0$. This is true for any $(k,l)\in\Lambda$, completing the proof of the Lemma.\qed
\vspace{0.1in}\\
We define a linear functional $\tau$ on $\clv_{2,+}$ by defining it on $\clb$ by the following prescription:\\
 $\tau(S_{e_{i}}S_{e_{j}}^{\ast})=\delta_{ij}$ for $(i,j)\in\Lambda$ and $\tau(p_{v_{k}})=1$ for all $k\in\cli$.
 It is clear from the definition that $\tau(S_{e_{i}}S_{e_{j}}^{\ast})=\delta_{ij}$ for all $i,j=1,\ldots,n$.
\blmma
For a connected graph $\Gamma=\{E=(e_{1},\ldots,e_{n}),V=(v_{1},\ldots,v_{m})\}$, let $F^{\Gamma}$ be an $n\times n$ matrix whose $ij$-th entries for all $i,j=1,\ldots,n$ are given by $\tau(S_{e_{i}}^{\ast}S_{e_{j}})$. Then $F^{\Gamma}$ is an invertible matrix.
\elmma
{\bf Proof}:\\
Since $S_{e_{i}}^{\ast}S_{e_{j}}=0$ for $i\neq j$, $F^{\Gamma}$ is a diagonal matrix. So in order to show that $F^{\Gamma}$ is invertible, it is enough to show that the diagonal entries are all non zero i.e. $\tau(S_{e_{i}}^{\ast}S_{e_{i}})\neq 0$ for all $i=1,\ldots,n$. For a fixed $i$, $s^{-1}(t(e_{i}))=\emptyset$ or $s^{-1}(t(e_{i}))\neq\emptyset$. If $s^{-1}(t(e_{i}))=\emptyset$, then $t(e_{i})\in\cli$ and hence by definition of $\tau$, $\tau (S_{e_{i}}^{\ast}S_{e_{i}})=\tau(p_{t(e_{i})})=1$. If $s^{-1}(t(e_{i}))\neq\emptyset$, then we assume that there are $l$($l\geq 1$) number of edges $\{e_{k}:k=1,\ldots,l\}$ such that $t(e_{i})=s(e_{k})$ for all $k=1,\ldots,l$. Then by the defining relations of the graph $C^{\ast}$-algebra, we have $S_{e_{i}}^{\ast}S_{e_{i}}=\sum_{k=1}^{l}S_{e_{k}}S_{e_{k}}^{\ast}$. Hence $\tau(S_{e_{i}}^{\ast}S_{e_{i}})=l$. In any case $\tau(S_{e_{i}}^{\ast}S_{e_{i}})\neq 0$, completing the proof of the Lemma.\qed

\subsection{Quantum symmetry}
\bdfn
We call a faithful action $\alpha$ of a CQG $\clq$ on a connected graph $C^{\ast}$-algebra $C^{\ast}(\Gamma)$ corresponding to a graph $\Gamma=\{E=(e_{1},\ldots,e_{n}),V=(v_{1},\ldots,v_{m})\}$ linear if $\alpha(S_{e_{i}})=\sum_{j=1}^{n}S_{e_{j}}\ot q_{ji}$, where $q_{ij}\in\clq$ for $i,j=1,...,n$.
\edfn
\brmrk
Note that since the graphs are connected, it suffices to define actions on the partial isometries corresponding to the edges in order to get the action on the generators of the graph $C^{\ast}$-algebra.
\ermrk
Recall the vector space $\clv_{2,+}$ from the previous subsection. 
\blmma
A linear action $\alpha$ of a CQG $(\clq,\Delta)$ on $C^{\ast}(\Gamma)$ preserves the vector space $\clv_{2,+}$ i.e. $\alpha(\clv_{2,+})\subset \clv_{2,+}\ot\clq$.
\elmma
{\bf Proof}:\\
Note that it suffices to show that $\alpha$ maps the basis elements $\clv_{2,+}$ into $\clv_{2,+}\ot\clq$. For that, writing $\alpha(S_{e_{i}})=\sum_{j}{S_{e_{j}}}\ot q_{ji}$, we see that $\alpha(S_{e_{i}}S_{e_{j}}^{\ast})=\sum_{k,l}S_{e_{k}}S_{e_{l}}^{\ast}\ot q_{ki}q_{lj}^{\ast}\in\clv_{2,+}\ot\clq$. For some vertex $p_{v_{l}}$ which is not a source of any edge, by assumption there exists some $e_{k}$ such that $S_{e_{k}}^{\ast}S_{e_{k}}=p_{v_{l}}$. Then $\alpha(p_{v_{l}})=\sum_{j}S_{e_{j}}^{\ast}{S_{e_{j}}}\ot q_{jk}^{\ast}q_{jk}$. Now using the fact that $S_{e_{i}}^{\ast}S_{e_{i}}=\sum_{s(e_{l})=t(e_{i})}S_{e_{l}}S_{e_{l}}^{\ast}$, we finish the proof.\qed 
\bdfn
For a finite graph $\Gamma$ we define a category $\clc^{\rm Lin}_{\tau}$ whose objects are $((\clq,\Delta),\alpha)$, where $\alpha$ is a {\bf linear} faithful $C^{\ast}$-action of a CQG $(\clq,\Delta)$ on $C^{\ast}(\Gamma)$ such that $\alpha$ preserves $\tau$ on $\clv_{2,+}$. Morphism between two objects $((\clq_{1},\Delta_{1}),\alpha_{1})$ and $((\clq_{2},\Delta_{2}),\alpha_{2})$ in the category $\clc^{\rm Lin}_{\tau}$ is given by a CQG morphism $\Phi:\clq_{1}\raro\clq_{2}$ such that $({\rm id}\ot\Phi)\alpha_{1}=\alpha_{2}$. 
\edfn
For a connected graph $\Gamma$, recall the invertible, diagonal matrix $F^{\Gamma}$ from Subsection \ref{graph}. We shall denote the $ij$-th entry of $F^{\Gamma}$ by $f_{ij}$.
\bppsn
\label{U_n}
For a finite, connected graph $\Gamma$, there is a surjective $C^{\ast}$-homomorphism from $A_{u^{t}}(F^{\Gamma})$ to any object in the category $\clc^{\rm Lin}_{\tau}$. 
\eppsn
{\bf Proof}:\\
As before we fix a finite, connected graph $\Gamma=\{E=(e_{1},...,e_{n}),V=(v_{1},...,v_{m})\}$. Let $\alpha$ be a linear action of a CQG $(\clq,\Delta)$ in the category $\clc^{\rm Lin}_{\tau}$ on $C^{\ast}(\Gamma)$ given by
\begin{displaymath}
\alpha(S_{e_{i}})=\sum_{k=1}^{n}S_{e_{k}}\ot u_{ki},
\end{displaymath}
for $i=1,...,n$. Using the fact that $\alpha$ preserves $\tau$ on $\clv_{2,+}$, we have for all $i,j=1,...,n$,
\begin{eqnarray}
&&(\tau\ot {\rm id})\alpha(S_{e_{i}}S_{e_{j}}^{\ast})=\tau(S_{e_{i}}S_{e_{j}}^{\ast})1 \nonumber\\ 
&\Rightarrow& \sum_{k,l=1}^{n}\tau(S_{e_{k}}S_{e_{l}}^{\ast})u_{ki}u_{lj}^{\ast}=\delta_{ij}1 \nonumber\\
\label{5} &\Rightarrow& \sum_{k=1}^{n}u_{ki}u_{kj}^{\ast}=\delta_{ij}1 
\end{eqnarray} 
If we denote the matrix $((u_{ij}))_{i,j=1,\ldots,n}$ by $Q$, then equation (\ref{5}) implies that $Q^{t}Q^{t^{\ast}}={\rm Id}_{n\times n}$. Since we are considering CQG action, the invertibility of $Q^{t}$ ensures that $Q^{t}$ is unitary. Again using the $\tau$-preserving property on $\clv_{2,+}$, we get
\begin{eqnarray}
&& (\tau\ot {\rm id})\alpha(S_{e_{i}}^{\ast}S_{e_{j}})=\tau(S_{e_{i}}^{\ast}S_{e_{j}})1 \nonumber\\
&\Rightarrow& \sum_{k=1}^{n}\tau(S_{e_{k}}^{\ast}S_{e_{k}})u_{ki}^{\ast}u_{kj}=f_{ij}\nonumber\\
&\Rightarrow& \label{6} \sum_{k=1}^{n} u_{ki}^{\ast}f_{kk}u_{kj}=f_{ij}
\end{eqnarray}
Equation (\ref{6}) can be written in the matrix equation form by $Q^{\ast}F^{\Gamma}Q=F^{\Gamma}$. Hence using the invertibilty of the matrix $Q$, we conclude that $Q^{-1}=F^{{\Gamma}^{-1}}Q^{\ast}F^{\Gamma}$. Then the claim of the proposition follows from the universal property of $A_{u^{t}}(F^{\Gamma})$.\qed\vspace{0.1in}\\
Adapting the arguments as in Theorem 4.8 of \cite{univ}, we prove the following 
\bthm
\label{univ}
For a finite, connected graph $\Gamma$, the category $\clc^{\rm Lin}_{\tau}$ admits a universal object.
\ethm
{\bf Proof}:\\
For a finite, connected graph $\Gamma=\{E=(e_{1},...,e_{n}),V=(v_{1},...,v_{m})\}$, let us define a category $\clc$ with objects $(\alpha_{C},C)$, where $C$ is a $C^{\ast}$-algebra generated by $\{t_{ij}: i,j=1,...,n\}$ with both the $C$-valued matrix $T=((t_{ij}))$ and $T^{t}=((t_{ij}))^{t}$ satisfying the defining relations of $A_{u^{t}}(F^{\Gamma})$ and $\alpha_{C}:C^{\ast}(\Gamma)\raro C^{\ast}(\Gamma)\ot C$ is a $C^{\ast}$-homomorphism. By definition of the category $\clc$, it is clear that for any object $(\alpha_{C},C)\in \clc$, there is a surjective $C^{\ast}$-homomorphism from $A_{u^{t}}(F^{\Gamma})$ to $C$ sending $u_{ij}$ to $t_{ij}$. Let $I_{C}$ be the kernel of the $C^{\ast}$-homomorphism and $I=\{\cap_{C} I_{C}:(\alpha_{C},C)\in\clc\}$. Clearly $I$ is a closed $C^{\ast}$-ideal of $A_{u^{t}}(F^{\Gamma})$. We define $\clu=A_{u^{t}}(F^{\Gamma})/I$.  We denote the canonical projection from $A_{u^{t}}(F^{\Gamma})$ onto $\clu$ by $P$ and $P(u_{ij})$ by $\overline{u_{ij}}$. We define a map $\alpha:C^{\ast}(\Gamma)\raro C^{\ast}(\Gamma)\ot \clu$ by $\alpha(S_{e_{i}})=\sum_{j=1}^{n}S_{e_{j}}\ot\overline{u_{ij}}$. Then $\alpha$ is a $C^{\ast}$-homomorphism and then it is clear from the definition that $\clu$ is the universal object in the category $\clc$. In the following we shall equip $\clu$ with a CQG structure and subsequently prove that $\clu$ with the given structure is in deed the universal object in the category $\clc^{\rm Lin}_{\tau}$.\\
\indent It is clear that it suffices to define coproduct on $\overline{u_{ij}}$ for $i,j=1,...,n$. To that end we define $\overline{U_{ij}}=\sum_{k=1}^{n}\overline{u_{ik}}\ot \overline{u_{kj}}\in \clu\ot \clu$ for $i,j=1,...,n$. Clearly $\overline{U_{ij}}=(P\ot P)\Delta_{0}(u_{ij})$, where $\Delta_{0}$ is the coproduct of $A_{u^{t}}(F^{\Gamma})$. Now we define $\beta: C^{\ast}(\Gamma)\raro C^{\ast}(\Gamma)\ot \clu\ot \clu$ by $\beta(S_{e_{i}})=\sum_{j=1}^{n}S_{e_{j}}\ot\overline{U_{ji}}$. It is easy to see that $\beta=(\alpha\ot {\rm id})\alpha$ proving that $\beta$ is a $C^{\ast}$-homomorphism. Also $((\overline{U_{ij}}))$ and $((\overline{U_{ij}}))^{t}$ satisfy the defining relations of $A_{u^{t}}(F^{\Gamma})$ proving that $\clu\ot \clu$ is also an object in the category $\clc$. Hence by universal property of $\clu$, there is a surjective $C^{\ast}$-homomorphism $\Delta:\clu\raro \clu\ot \clu$ sending $\overline{u_{ij}}$ to $\overline{U_{ij}}=\sum_{k=1}^{n}\overline{u_{ik}}\ot\overline{u_{kj}}$. We claim that $\Delta$ is the required coproduct on $\clu$. It is straightforward to check that $(\Delta\ot{\rm id})\Delta=({\rm id}\ot\Delta)\Delta$. Using the condition $\Delta \circ P= (P \ot P)\circ \Delta_{0}$ we  can conclude that $\Delta_0$ maps $Ker(P)$ to $Ker(P \ot P)$. Moreover, from the density of each of the linear spans of $\Delta_{0}(\clq_0)(1\ot \clq_0)$ as well as $\Delta_0(\clq_0)(\clq_0 \ot 1)$ in $A_{u^{t}}(F^{\Gamma}) \ot A_{u^{t}}(F^{\Gamma})$, where $\clq_0$ is the *-algebra generated by $u_{ij}$'s, we get (applying $P \ot P$) similar density with $A_{u^{t}}(F^{\Gamma})$ and $\Delta_0$ replaced by $\clu$ and $\Delta$ respectively, and $\clq_0$ by the algebra generated by $\overline{u_{ij}}$'s. Therefore, $Ker(P)$ is a closed Hopf ideal and $\clu$ becomes a quantum subgroup of $A_{u^{t}}(F^{\Gamma})$ with coproduct $\Delta$.\\
\indent Now we finish the proof of the theorem by proving that $((\clu,\Delta),\alpha)$ is the universal object in the category $\clc^{\rm Lin}_{\tau}$. $\alpha$ is a linear faithful $C^{\ast}$-action of $(\clu,\Delta)$ on $C^{\ast}(\Gamma)$ by definition of $\alpha$ and $\Delta$. Also using the fact that $((\overline{U_{ij}}))$ and $((\overline{U_{ij}}))^{t}$ satisfy the defining relation of $A_{u^{t}}(F^{\Gamma})$ it is easy to see that $\alpha$ preserves $\tau$ on $\clv_{2,+}$. Hence $((\clu,\Delta),\alpha)$ in deed belongs to the category $\clc^{\rm Lin}_{\tau}$. Now given any object $((Q,\Delta_{Q}),\alpha_{Q})$ in the category $\clc^{\rm Lin}_{\tau}$ such that $\alpha_{Q}(S_{e_{i}})=\sum_{j=1}^{n}S_{e_{j}}\ot q_{ji}$, it follows from Proposition \ref{U_n}, $(\alpha_{Q},Q)$ is an object in the category $\clc$ (forgetting the CQG structure of $(Q,\Delta_{Q})$). So by universality of $\clu$, there is a surjective $C^{\ast}$-homomorphism $\Phi:\clu\raro Q$ sending $\overline{u_{ij}}$ to $q_{ij}$. It is clear from the fact $\Delta_{Q}(q_{ij})=\sum_{k=1}^{n}q_{ik}\ot q_{kj}$ that in fact $\Phi$ is a CQG morphism. Also by definition $({\rm id}\ot \Phi)\alpha=\alpha_{Q}$ which proves that $\clu$ is the universal object in the category $C^{\rm Lin}_{\tau}$. This completes the proof of the theorem. \qed 
\vspace{0.1in}\\
From now on, we denote the universal object in the category $\clc^{\rm Lin}_{\tau}$ by $Q^{\rm Lin}_{\tau}$.   
\subsection{Non rigidity of quantum symmetry}
In this subsection we shall show that the quantum symmetry of a finite, connected graph $\Gamma=\{E=(e_{1},\ldots,n),V=(v_{1},\ldots, v_{m})\}$ in our sense is quite large. First we show that if the underlying graph has no multiple edges or loops, then the quantum symmetry of the underlying directed graph in the sense of Banica belongs to the category $\clc^{\rm Lin}_{\tau}$. We denote the vertex by the index set itself for simplicity, i.e. we denote the vertex $v_{i}$ by $i$ for $i=1,...,m$. Also the corresponding projections will be denoted by $p_{i}$. Recall $Q^{aut}_{Ban}(\Gamma)$ from Subsection \ref{graph_bich}.
\bthm
\label{ban}
There is a $C^{\ast}$-action $\alpha$ of $(Q^{aut}_{Ban}(\Gamma),\Delta)$ on $C^{\ast}(\Gamma)$ such that $((Q^{aut}_{Ban}(\Gamma),\Delta),\alpha)$ is an object in the category $\clc^{\rm Lin}_{\tau}$.
\ethm
{\bf Proof}:\\
We define a map $\alpha:C^{\ast}(\Gamma)\raro C^{\ast}(\Gamma)\ot Q^{aut}_{Ban}(\Gamma)$ by the following
\begin{eqnarray*}
&&\alpha(p_{i})=\sum_{k=1}^{m}p_{k}\ot u_{ik}\\
&& \alpha(S_{e_{j}})=\sum_{l=1}^{n}S_{e_{l}}\ot u_{s(e_{j})s(e_{l})}u_{t(e_{j})t(e_{l})}. 
\end{eqnarray*}
 By definition $\alpha$ is linear. So in order to prove the theorem, we have to show that $\alpha$ is a well defined faithful $C^{\ast}$-action and it preserves $\tau$ on $\clv_{2,+}$. For the proof that $\alpha$ is a faithful $C^{\ast}$-action we refer the reader to Theorem 4.1 of \cite{Web}.
Now we proceed to prove that $\alpha$ is $\tau$-preserving on $\clv_{2,+}$. For $i=1,...,n$,
 \begin{eqnarray*}
(\tau\ot {\rm id})\alpha(S_{e_{i}}S_{e_{i}}^{\ast})&=&\sum_{l,k=1}^{n}\tau(S_{e_{l}}S_{e_{k}}^{\ast})u_{s(e_{i})s(e_{l})}u_{t(e_{i})t(e_{l})}u_{t(e_{i})t(e_{k})}u_{s(e_{i})s(e_{k})}\\
&=&\sum_{k=1}^{n} u_{s(e_{i})s(e_{k})}u_{t(e_{i})t(e_{k})}u_{t(e_{i})t(e_{k})}u_{s(e_{i})s(e_{k})}\\
&=& \sum_{k=1}^{n} u_{s(e_{i})s(e_{k})}u_{t(e_{i})t(e_{k})}u_{s(e_{i})s(e_{k})} \ (using \ u_{ij}^2=u_{ij})\\
\end{eqnarray*}
Now using the fact that for $(i^{\prime},j^{\prime})\not\in E$, $u_{s(e_{i})i^{\prime}}u_{t(e_{i})j^{\prime}}u_{s(e_{i})i^{\prime}}=0$ (equation (\ref{2})), we get the last expression equals to $\sum_{i^{\prime},j^{\prime}}u_{s(e_{i})i^{\prime}}u_{t(e_{i})j^{\prime}}u_{s(e_{i})i^{\prime}}$. Using equation (\ref{1}), we get this expression is equal to $1$. Hence $(\tau\ot {\rm id})\alpha(S_{e_{i}}S_{e_{i}}^{\ast})=\tau(S_{e_{i}}S_{e_{i}}^{\ast})$ for all $i=1,\ldots,n$. \\
\indent For $i\neq j$, if $S_{e_{i}}S_{e_{j}}^{\ast}=0$, then $\alpha$ is automatically $\tau$ preserving for such elements. So for $i\neq j$, we assume $S_{e_{i}}S_{e_{j}}^{\ast}\neq 0$ i.e. $t(e_{i})=t(e_{j})$. For such elements using the fact that $u_{t(e_{i})t(e_{k})}=u_{t(e_{j})t(e_{k})}$, we have 
\begin{eqnarray*}
(\tau\ot{\rm id})\alpha(S_{e_{i}}S_{e_{j}}^{\ast})&=&\sum_{k=1}^{n}u_{s(e_{i})s(e_{k})}u_{t(e_{i})t(e_{k})}u_{t(e_{j})t(e_{k})}u_{s(e_{j})s(e_{k})}\\
&=& \sum_{k=1}^{n} u_{s(e_{i})s(e_{k})}u_{t(e_{i})t(e_{k})}u_{s(e_{j})s(e_{k})}\\
\end{eqnarray*}
Using similar arguments as before, we can show that the last expression equals to $\sum_{i^{\prime}}u_{s(e_{i})i^{\prime}}u_{s(e_{j})i^{\prime}}$. Since $\Gamma$ does not have any multiple edge, $t(e_{i})=t(e_{j})\Rightarrow s(e_{i})\neq s(e_{j})$. So using equation (\ref{1}), we have for all $i^{\prime}$, $u_{s(e_{i})i^{\prime}}u_{s(e_{j})i^{\prime}}=0$. Hence we have for all $i,j=1,...,n$,
\begin{displaymath}
(\tau\ot{\rm id})\alpha(S_{e_{i}}S_{e_{j}}^{\ast})=\tau(S_{e_{i}}S_{e_{j}}^{\ast})1.
\end{displaymath}
To complete the proof of $\tau$-preserving property we need to show that $(\tau\ot{\rm id})(\alpha(p_{i}))=\tau(p_{i}).1$ for those $i$'s which are not source of any vertex. Hence we assume $p_{i}=S_{e_{k}}^{\ast}S_{e_{k}}$ for some edge $e_{k}$. By definition of $\tau$, $\tau(p_{i})=1$. So we need to show that $(\tau\ot{\rm id})\alpha(S_{e_{k}}^{\ast}S_{e_{k}})=1$. 
\begin{eqnarray*}
(\tau\ot{\rm id})\alpha(S_{e_{k}}^{\ast}S_{e_{k}})&=&\sum_{l=1}^{n}\tau(S_{e_{l}}^{\ast}S_{e_{l}}) u_{t(e_{k})t(e_{l})}u_{s(e_{k})s(e_{l})} u_{t(e_{k})t(e_{l})}\\
\end{eqnarray*}
Let $e_{l}\in E$ be such that $s^{-1}(t(e_{l}))\neq \emptyset$. Then there is some $f\in E$ such that $t(e_{l})=s(f)$. Hence for such $l$, $u_{t(e_{k})t(e_{l})}=u_{is(f)}$. But since $s^{-1}(i)=\emptyset$, by Lemma \ref{Weber}, $u_{is(f)}=0$ i.e. for $e_{k}, e_{l}\in E$ such that $s^{-1}(t(e_{k}))=\emptyset$ and $s^{-1}(t(e_{l}))\neq \emptyset$, $u_{t(e_{k})t(e_{l})}=0$. Hence using the fact that $\tau(S_{e_{l}}^{\ast}S_{e_{l}})=1$ for all $e_{l}$ with $s^{-1}(t(e_{l}))=\emptyset$, we get 
\begin{eqnarray*}
 && (\tau\ot{\rm id})\alpha(S_{e_{k}}^{\ast}S_{e_{k}})\\
 &=&\sum_{s^{-1}(t(e_{l}))=\emptyset}\tau(S_{e_{l}}^{\ast}S_{e_{l}})u_{t(e_{k})t(e_{l})}u_{s(e_{k})s(e_{l})}u_{t(e_{k})t(e_{l})}\\
 &=& \sum_{s^{-1}(t(e_{l}))=\emptyset} u_{t(e_{k})t(e_{l})}u_{s(e_{k})s(e_{l})}u_{t(e_{k})t(e_{l})}.
 \end{eqnarray*}
 Again using similar arguments the last expression equals to $\sum_{k^{\prime},j^{\prime}}u_{t(e_{k})k^{\prime}}u_{s(e_{k})j^{\prime}}u_{t(e_{k})k^{\prime}}=1$, proving the $\tau$-preserving property and hence completing the proof of the Theorem.\qed
 \vspace{0.1in}\\
Now we show that in fact almost always we have genuine quantum symmetry of a finite, connected graph. Note that for the following Proposition, the condition on the graph having no {\bf multiple edge} or {\bf loop} is relaxed. 
\bppsn
\label{non_rigid}
For a finite, connected graph $\Gamma$ with $n$ edges, $(\underbrace{(C(S^{1})\ast\ldots\ast C(S^{1})}_{n-copies},\Delta),\alpha)$ always belongs to the category $\clc^{\rm Lin}_{\tau}$ where $\alpha(S_{e_{i}})=S_{e_{i}}\ot q_{i}$ for all $i=1,\ldots,n$ and $q_{i}$'s are unitaries generating the free product.
\eppsn
{\bf Proof}:\\
We shall first prove the existence of such a $C^{\ast}$-homomorphism $\alpha$ by proving that $S_{e_{i}}^{\prime}:=S_{e_{i}}\ot q_{i}$'s also satisfy the same relations of $C^{\ast}(\Gamma)$. Then by universal property of $C^{\ast}(\Gamma)$, we can establish the existence of such a well defined $C^{\ast}$-homomorphism. We extend the definition of $\alpha$ on the vertex set by $\alpha(p_{i})=p_{i}\ot 1$ for $i=1,\ldots,m$. We denote $p_{i}\ot 1$ by $p_{i}^{\prime}$. Now we show that $S_{e_{j}}^{\prime},p_{i}^{\prime}$ for $i=1,\ldots,m$ and $j=1,\ldots,n$ satisfies the relations in $C^{\ast}(\Gamma)$. Firstly observe that $p_{i}^{\prime}p_{j}^{\prime}=\delta_{ij}p_{i}^{\prime}$ and $p_{i}^{\prime\ast}=p_{i}^{\prime}$ and $\sum_{i=1}^{m}p_{i}^{\prime}=1$. Also using the unitarities of $q_{i}$'s, we have $S_{e_{i}}^{\prime\ast}S_{e_{i}}^{\prime}=S_{e_{i}}^{\ast}S_{e_{i}}\ot 1=p_{t(e_{i})}\ot 1=p_{t(e_{i})}^{\prime}$ and $\sum_{s(e_{i})=l}S_{e_{i}}^{\prime}S_{e_{i}}^{\prime\ast}=\sum_{s(e_{i})=l}S_{e_{i}}S_{e_{i}}^{\ast}\ot 1=p_{l}^{\prime}$, proving the existence of such a $C^{\ast}$-homomorphism $\alpha$. Coassociativity and faithfulness follow immediately. For the span density condition observe that $\alpha(S_{e_{i}})(1\ot q_{i}^{\ast})=S_{e_{i}}\ot 1$. So $(S_{e_{i}}\ot 1) \in \ {\rm Sp} \ \alpha(C^{\ast}(\Gamma))(1\ot \clq)$. Similarly $(S_{e_{i}}^{\ast}\ot 1)$ and $(p_{j}\ot 1)$ all belong to ${\rm Sp} \ \alpha(C^{\ast}(\Gamma))(1\ot \clq)$. Now the action being unital, the span density condition follows from standard facts of CQG actions. It is easy to see the $\tau$-preserving property, completing the proof of the Proposition.\qed
\vspace{0.1in}\\
The previous Proposition yields the following
\bcrlre
For $n\geq 2$, $Q^{\rm Lin}_{\tau}$ for a connected graph $\Gamma$ with $n$ number of edges is always non commutative as a $C^{\ast}$-algebra.
\ecrlre

\section{Computation of $Q^{\rm Lin}_{\tau}$ for some connected graphs} 
1. 
\begin{center} \begin{tikzpicture}
\begin{scope}[thick, every node/.style={sloped,allow upside down}]
\draw (0,7.5) -- node {\midarrow} (2,7.5);
\draw (2,7.5) -- node {\midarrow} (4,7.5);
\filldraw [black] (0,7.5) circle (1pt);
\filldraw [black] (2,7.5) circle (1pt);
\filldraw [black] (4,7.5) circle (1pt);
\filldraw [black] (4.2,7.5) circle (.25pt);
\filldraw [black] (4.4,7.5) circle (.25pt);
\filldraw [black] (4.6,7.5) circle (.25pt);
\filldraw [black] (4.8,7.5) circle (.25pt);
\filldraw [black] (5,7.5) circle (1pt);
\draw (5,7.5) -- node {\midarrow} (7,7.5);
\end{scope}
\end{tikzpicture}
\end{center}
\indent Let us consider the above graph $\Gamma$ with $n$ consecutive edges and $(n+1)$ vertex joining them. We denote the projections corresponding to the vertices by $\{p_{i}\}_{i=0,\ldots,n}$ and the partial isometries corresponding to the edges by $\{S_{i}\}_{i=1,\ldots,n}$. We have the following relations in the graph $\Gamma$:
\begin{eqnarray}
 && S_{1}S_{1}^{\ast}=p_{0}, S_{n}^{\ast}S_{n}=p_{n},\sum_{i=0}^{n}p_{i}=1 \label{8}\\
 && S_{i}^{\ast}S_{i}=S_{i+1}S_{i+1}^{\ast}=p_{i}, i=1,\ldots,(n-1) \label{9}\\
 && S_{i}^{\ast}S_{j}=S_{i}S_{j}^{\ast}=0, i,j=1,...,n \ and \ i\neq j \label{10}.
\end{eqnarray}
\bthm
 For the graph $\Gamma$, $Q^{\rm Lin}_{\tau}$ is isomorphic to $(\underbrace{C(S^{1})\ast\ldots\ast C(S^{1})}_{n-times},\Delta)$.  
\ethm
We start with the following
\blmma
Let $((\clq,\Delta),\alpha)$ be an object in the category $\clc^{\rm Lin}_{\tau}$ for the graph $\Gamma$ so that $\alpha(S_{i})=\sum_{k=1}^{n}S_{k}\ot q_{ki}$ for $i=1,\ldots,n$ and $((q_{ki}))_{i,k=1,...,n}\in M_{n}(\clq)$. 
Then $q_{ij}=0$ for $i\neq j$.
\elmma
{\bf Proof}:\\
We have
\begin{eqnarray*}
 \alpha(S_{i}^{\ast}S_{i})=\sum_{k=1}^{n}S_{k}^{\ast}S_{k}\ot q_{ki}^{\ast}q_{ki}\\
 \alpha(S_{i}S_{i}^{\ast})=\sum_{k=1}^{n}S_{k}S_{k}^{\ast}\ot q_{ki}q_{ki}^{\ast}.
\end{eqnarray*}
Using the relation $\alpha(S_{i}^{\ast}S_{i})=\alpha(S_{i+1}S_{i+1}^{\ast})$ for $i=1,\ldots,(n-1)$, we get
\begin{eqnarray*}
 &&\sum_{k=1}^{n}S_{k}^{\ast}S_{k}\ot q_{ki}^{\ast}q_{ki}=\sum_{k=1}^{n}S_{k}S_{k}^{\ast}\ot q_{ki+1}q_{ki+1}^{\ast}\\
 &\Rightarrow& S_{1}S_{1}^{\ast}\ot q_{1i+1}q_{1i+1}^{\ast}-S_{n}^{\ast}S_{n}\ot q_{ni}^{\ast}q_{ni}+\sum_{k=1}^{n-1}S_{k}^{\ast}S_{k}\ot(q_{k+1i+1}q_{k+1i+1}^{\ast}-q_{ki}^{\ast}q_{ki})=0
\end{eqnarray*}
Using the orthogonality of the projections $\{p_{i}\}_{i=0,\ldots,n}$, we have for all $i,k=1,\ldots,(n-1)$
\begin{eqnarray}
 && q_{1i+1}=q_{ni}=0 \label{11}\\
 && q_{k+1i+1}q_{k+1i+1}^{\ast}=q_{ki}^{\ast}q_{ki} \label{12}
\end{eqnarray}
By (\ref{11}), $q_{1i}=0$ for all $i=2,\ldots,n$. Using this and (\ref{12}), recursively, we get $q_{li}=0$ for all $l=1,\ldots,(n-1)$ and $i=l+1,\ldots,n$. Now applying the antipode to $q_{li}$ for all $l=1,\ldots,n$ and $i=l+1,\ldots,n$, we show that $q_{ij}=0$ for $i\neq j$. \qed
\vspace{0.1in}\\
{\bf Proof of the Theorem}:\\
Note that by Proposition \ref{non_rigid}, it is enough to show that we have a CQG homomorphism from $(\underbrace{C(S^{1})\ast\ldots\ast C(S^{1})}_{n-times},\Delta)$ onto $Q^{\rm Lin}_{\tau}$. If $\alpha$ is the action of $Q^{\rm Lin}_{\tau}$ given by $\alpha(S_{e_{i}})=\sum_{j}S_{e_{j}}\ot q_{ji}$, then by the previous Lemma we have $\alpha(S_{i})=S_{i}\ot q_{ii}$ for $q_{ii}\in Q^{\rm Lin}_{\tau}$ for $i=1,\ldots,n$. From the relation (\ref{8}), (\ref{9}) it follows that $q_{i+1i+1}q_{i+1i+1}^{\ast}=q_{ii}^{\ast}q_{ii}$ for all $i=1,\ldots,n$. Using the fact that $\sum_{i=0}^{n}p_{i}={\rm id}$, we get the following
\begin{eqnarray}
 && \sum_{i=1}^{n}\alpha(S_{i}S_{i}^{\ast})+\alpha(S_{n}^{\ast}S_{n})=1\ot 1 \nonumber \\
 &\Rightarrow& \sum_{i=1}^{n}S_{i}S_{i}^{\ast}\ot q_{ii}q_{ii}^{\ast}+ S_{n}^{\ast}S_{n}\ot q_{nn}^{\ast}q_{nn}=1\ot 1 \label{13}
 \end{eqnarray}
 Multiplying the last equation by $S_{k}$ from right and using the relation $S_{n}S_{k}=0$ for $k=1,\ldots,n$ and $S_{i}S_{i}^{\ast}S_{k}=\delta_{ik}S_{k}$, we get
 \begin{displaymath}
  S_{k}\ot q_{kk}q_{kk}^{\ast}=S_{k}\ot 1,
 \end{displaymath}
which implies that $q_{kk}q_{kk}^{\ast}=1$ for all $k=1,\ldots,n$. Using $q_{k+1k+1}q_{k+1k+1}^{\ast}=q_{kk}^{\ast}q_{kk}$ for $k=1,\ldots,(n-1)$, we get $q_{kk}^{\ast}q_{kk}=1$ for $k=1,\ldots,(n-1)$. Again multiplying equation (\ref{13}) by $S_{n}$ from left and using $S_{n}S_{i}=0$ for all $i$ we get $q_{nn}^{\ast}q_{nn}=1$ implying that all the $q_{ii}$'s are unitary. The action automatically preserves $\tau$. Hence by the universal property of $(\underbrace{C(S^{1})\ast\ldots\ast C(S^{1})}_{n-times},\Delta)$, we have a surjective CQG morphism from $(\underbrace{C(S^{1})\ast\ldots\ast C(S^{1})}_{n-times},\Delta)$ onto $Q^{\rm Lin}_{\tau}$ sending $z_{i}$ to $q_{ii}$ for $i=1,\ldots,n$, where $z_{i}$'s are the canonical generators of the free product. This completes the proof of the theorem. \qed
\brmrk
We observed in the previous Lemma, for the graph $\Gamma$, the presence of $\tau$ in the category $\clc^{\rm lin}_{\tau}$ is superfluous.
\ermrk
\vspace{0.25in}
2.  
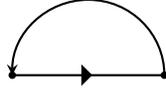
\begin{center} \begin{tikzpicture}
\begin{scope}[thick, every node/.style={sloped,allow upside down}]
\draw (0,7.5) -- node {\midarrow} (2,7.5);  
\filldraw [black] (0,7.5) circle (1pt);
\filldraw [black] (2,7.5) circle (1pt);
\draw [black, -stealth ] (2,7.5) arc (0:180:1) node [right] {};
\end{scope}
\end{tikzpicture}
\end{center}
\indent The above graph is a complete graph $\Gamma$ with two vertices $\{0,1\}$ and two edges $e_{1},e_{2}$. If we denote the projections corresponding to the vertices by $p_{0},p_{1}$ and the partial isometries corresponding to the edges by $S_{1},S_{2}$, then $C^{\ast}(\Gamma)$ is generated by $S_{1},S_{2}$ satisfying the following relations:
\begin{eqnarray}
 S_{1}^{\ast}S_{1}=S_{2}S_{2}^{\ast}=p_{1},S_{2}^{\ast}S_{2}=S_{1}S_{1}^{\ast}=p_{0}, 
\end{eqnarray}
where $p_{0},p_{1}$ are mutually orthogonal projections with $p_{0}+p_{1}=1$. Clearly this is an {\bf infinite} dimensional $C^{\ast}$-algebra. Let $((\clq,\Delta),\alpha)$ be an object in $\clc^{\rm Lin}_{\tau}$ such that for $((q_{ij}))_{i,j=1,2}\in M_{2}(\clq)$,
\begin{eqnarray}
 && \alpha(S_{i})=\sum_{j=1}^{2}S_{j}\ot q_{ji}, \ i=1,2.
\end{eqnarray}
\bthm\label{doubling_result}
For a complete graph $\Gamma$ with 2 vertices, $Q^{\rm Lin}_{\tau}$ is isomorphic to $(\cld_{\theta}(C(S^{1})\ast C(S^{1})),\tilde{\Delta})$ (see Subsection \ref{qaut}) where $\theta$ is the order $2$ automorphism given by $\theta(z_{1})=z_{2},\theta(z_{2})=z_{1}$ with $z_{1},z_{2}$ being the canonical generators of $C(S^{1})\ast C(S^{1})$.
\ethm
We break down the proof in a few Lemmas.
\blmma
\label{alzero1}
\begin{eqnarray}
 && q_{11}^{\ast}q_{11}=q_{22}q_{22}^{\ast},q_{21}^{\ast}q_{21}=q_{12}^{\ast}q_{12} \nonumber\\
 && q_{11}q_{11}^{\ast}=q_{22}^{\ast}q_{22},q_{21}q_{21}^{\ast}=q_{12}^{\ast}q_{12}. \label{16}
\end{eqnarray}
\elmma
{\bf Proof}:\\
The orthogonality of $p_{0},p_{1}$ along with the relations $\alpha(S_{1}^{\ast}S_{1})=\alpha(S_{2}S_{2}^{\ast})$ and $\alpha(S_{2}^{\ast}S_{2})=\alpha(S_{1}S_{1}^{\ast})$ produce the relations (\ref{16}).\qed
\blmma
\label{alzero}
\begin{eqnarray}
 && q_{11}q_{12}=q_{12}q_{11}=q_{21}q_{11}=q_{11}q_{21}=0 \nonumber\\
 && q_{22}q_{12}=q_{12}q_{22}=q_{22}q_{21}=q_{21}q_{22}=0 
\end{eqnarray}
\elmma
{\bf Proof}:\\
Using the fact that $\alpha(S_{1}^{\ast}S_{2})=0$, we have the following relation
\begin{displaymath}
 p_{1}\ot q_{11}^{\ast}q_{12}+p_{0}\ot q_{21}^{\ast}q_{22}=0.
\end{displaymath}
Now using the orthogonality of $p_{0}$ and $p_{1}$, we can conclude that $q_{11}^{\ast}q_{12}=q_{21}^{\ast}q_{22}=0$.
\begin{eqnarray*}
 (q_{11}q_{21})^{\ast}(q_{11}q_{21})&=&q_{21}^{\ast}q_{11}^{\ast}q_{11}q_{21}\\
 &=& q_{21}^{\ast}q_{22}q_{22}^{\ast}q_{21} \ (by \ (\ref{16}))\\
 &=& 0 \ (q_{21}^{\ast}q_{22}=0)
\end{eqnarray*}
Hence $q_{11}q_{21}=0$. Applying $\kappa$, we get $q_{11}q_{12}=0$. With similar arguments we can show the rest of the relations.\qed
\blmma
\label{alzero2}
$q_{ij}$'s are normal and partial isometries for all $i,j=1,2$.
\elmma
{\bf Proof}:\\
Observe that the matrix $F^{\Gamma}$ is the identity matrix for this graph. So both $((q_{ij}))$ and $((q_{ij}))^{t}$ are unitaries implying the following:
\begin{eqnarray*}\label{eqn new}
  q_{11}q_{11}^{\ast}+q_{12}q_{12}^{\ast}= q_{11}^{\ast}q_{11}+q_{12}^{\ast}q_{12}=1.
\end{eqnarray*}
Using Lemma (\ref{alzero}), we get 
\begin{eqnarray}
 && q_{11}^{\ast 2}q_{11}=q_{11}^{\ast} \label{18}\\
 && q_{11}q_{11}^{\ast 2}=q_{11}^{\ast} \label{19}
\end{eqnarray}
Multiplying (\ref{18}) by $q_{11}$ from left and using (\ref{19}), we get $q_{11}^{\ast}q_{11}=q_{11}q_{11}^{\ast}$, proving that $q_{11}$ is normal. Similarly we can show that the rest of the elements are normal as well. Moreover, note that $q_{ij}$'s are partial isometries using the unitarity condition of $((q_{ij}))$ and $q_{ij}q_{ik}^*=0$ for all $i,j,k$ with $j\neq k$. \qed\\
Observe that all $q_{ij}q_{ij}^*$ are central projections. Now we turn to the proof of the theorem.
\vspace{0.1in}\\
{\bf Proof of the Theorem}:\\
Note that the underlying $C^*$-algebra of $Q_{\rm Lin}^{\tau}$ is generated by $\{q_{ij}\}_{i,j=1,2}$ satisfying all the relations of  Lemmas \ref{alzero}, \ref{alzero1}, \ref{alzero2} and $((q_{ij})), ((q_{ij}^t))$ are unitaries. Now we define a map $\phi:Q^{\rm Lin}_{\tau}\raro (C(S^{1})\ast C(S^{1}))\oplus (C(S^{1})\ast C(S^{1}))$ on the generators by $\phi(q_{11})=(z_{1},0), \phi(q_{12})=(0,z_{1}),\phi(q_{21})=(0,z_{2}),\phi(q_{22})=(z_{2},0)$, where $z_{1},z_{2}$ are canonical generators of the free product $C(S^{1})\ast C(S^{1})$. This map is clearly a $C^*$-isomorphism between them. Indeed, $Q_{\rm Lin}^{\tau}$  is isomorphic to $(\cld_{\theta}(C(S^{1})\ast C(S^{1})),\tilde{\Delta})$ as a CQG with respect to $\theta$ given in Theorem \ref{doubling_result}. \qed
\brmrk
 $Q^{aut}_{Ban}(\Gamma)=\mathbb{Z}_{2}$ which is strictly smaller than the quantum automorphism determined here.
\ermrk
3. {\bf Cuntz algebra}: Recall the definition of Cuntz algebra $O_{n}$ which is a unital $C^{\ast}$-algebra generated by $n$ isometries $\{S_{i}\}_{i=1,...,n}$ satisfying the following relations
\begin{displaymath}
S_{i}^{\ast}S_{i}=\sum_{k=1}^{n}S_{k}S_{k}^{\ast}=1, i=1,...,n.
\end{displaymath}
\bppsn
For Cuntz algebra $O_{n}$, $Q^{\rm Lin}_{\tau}$ is isomorphic to $U_{n}^{+}$.
\eppsn
{\bf Proof}:\\
First observe that from the definition of $\tau$ on $\clv_{2,+}$ and the defining relations of $O_{n}$, $\tau(S_{i}^{\ast}S_{i})=n$ for all $i=1,\ldots,n$. Hence (denoting the graph of Cuntz algebra by the same notation $O_{n}$) the matrix $F^{O_{n}}$ is the matrix $n{\rm Id}_{n\times n}$ implying that if we denote the matrix $((q_{ij}))$ by $U$, then the action $\alpha$ of $Q^{\rm Lin}_{\tau}$ given by $\alpha(S_{i})=\sum_{j=1}^{n}S_{j}\ot q_{ji}$ for $q_{ji}\in Q^{\rm Lin}_{\tau}$ preserves $\tau$ on $\clv_{2,+}$ if and only if $U^{\ast}U=U^{t}U^{t\ast}={\rm Id}_{n\times n}$. Then by the universal property of $U^{+}_{n}$, there is a surjective homomorphism from $U^{+}_{n}$onto $Q^{\rm Lin}_{\tau}$ sending $u_{ij}$ to $q_{ij}$. So to prove the theorem it would suffice to show that $U^{+}_{n}$ belongs to the category $\clc^{\rm Lin}_{\tau}$. To that end we define $\alpha:O_{n}\raro O_{n}\ot U_{n}^{+}$ by $\alpha(S_{i})=\sum_{j=1}^{n}S_{j}\ot u_{ji}$ for $u_{ij}\in U^{+}_{n}$. First we prove that in deed we have a well defined $C^{\ast}$-homomorphism given by $\alpha$. Like before we shall show that $\{\alpha(S_{i})\}_{i=1,\ldots,n}$ again satisfy the defining relations of $O_{n}$. $\alpha(S_{i}^{\ast}S_{i})=\sum_{k}S_{k}^{\ast}S_{k}\ot q_{ki}^{\ast}q_{ki}=1\ot 1$. On the other hand 
\begin{eqnarray*}
 \sum_{l=1}^{n}\alpha(S_{l}S_{l}^{\ast})&=&\sum_{k,m,l}S_{k}S_{m}^{\ast}\ot q_{kl}q_{ml}^{\ast}\\
 &=&\sum_{k,m}S_{k}S_{m}^{\ast}\ot(\sum_{l}q_{kl}q_{ml}^{\ast})\\
 &=&\sum_{k=1}^{n}S_{k}S_{k}^{\ast}\ot 1 \ ({\rm since} \ \sum_{l}q_{kl}q_{ml}^{\ast}=\delta_{km})\\
 &=& 1\ot 1.
\end{eqnarray*}
This proves the existence of $\alpha$. Coassociativity, faithfulness and $\tau$-preserving properties are easy to check. Hence $U_{n}^{+}$ is an object in $\clc^{\rm Lin}_{\tau}$ which completes the proof of the proposition.\qed

\vspace{0.12 in}
We would like to end the paper with the following concluding remarks:
\brmrk
1. In the context of noncommutative geometry recently semifinite spectral triples on Graph $C^{\ast}$-algebras have been constructed (see \cite{Pask}). So it is interesting to explore the quantum isometry of such semifinite spectral triples. We would like to see where our framework of quantum symmetry stands in this context.
\vspace{0.1in}\\
2. In \cite{Ortho}, universal object in the category of CQG's respecting certain orthogonal filtration is shown to exist for a large class of $C^{\ast}$-algebras including the Cuntz-Krieger algebra. Also in that paper for general graph $C^{\ast}$-algebras the authors alluded to a category similar to ours. Objects of their category need to preserve a faithful state on the whole $C^{\ast}$-algebra. Although our linear functional does not extend as a faithful state on the whole $C^{\ast}$-algebra, but we believe that a modified version would extend as a faithful state on the whole algebra. But even with that modified state, a priori their category is more restrictive and hence the quantum symmetry group is smaller than the quantum symmetry group in our sense. But sometimes they can coincide, like in the case of Cuntz algebra, we can normalize our linear functional and extending on the fixed point algebra, we get the canonical KMS state on the Cuntz algebra after composing with the canonical expectation. In that case their category and our category with normalized functional in deed coincide. But we believe that in general for graph $C^{\ast}$-algebras the approaches are different.  
\ermrk 
{\bf Acknowledgement}: The work started during the first author's visit to NISER Bhubaneshwar. The first author would like to thank NISER Bhubaneshwar as well as Sutanu Roy for warm hospitality during that visit. We would like to thank Malay Ranjan Biswal for his help in drawing figures in the paper. We would also like to thank Moritz Weber for valuable comments which improved Theorem \ref{ban}.  
 
Soumalya Joardar \\
Theoretical Science Unit,\\ 
JNCASR, Bangalore-560064, India\\ 
email: soumalya.j@gmail.com 
\vspace{0.1in}\\
Arnab Mandal\\
School Of Mathematical Sciences\\
NISER, HBNI,  Bhubaneswar,  Jatni-752050, India\\
email: arnabmaths@gmail.com

\end{document}